\documentclass[twocolumn]{sq}

\usepackage{latexsym,amssymb}
\usepackage{amsmath, amsbsy}
\usepackage{amsopn, amstext}

\newcommand{\comment}[1]{}
\newcommand{\lmd}{\lambda}
\newcommand{\Lmd}{\Lambda}

\newcommand{\R}{{\mathbb R}}  %ams bold
\newcommand{\rref}[1]{(\ref{#1})}
\newcommand{\norm}[1]{\left\Vert#1\right\Vert}
\newcommand{\abs}[1]{\left\vert #1 \right\vert}

\def\thefootnote{*}

\newcommand{\edo}{\end{document}}

\begin{document}

\begin{frontmatter}

\title{
\bf A Note on Overshoot Estimation in Pole Placements}
%\thanksref{footnoteinfo} }

\thanks[footnoteinfo]{This work was supported partly by the Chinese National Natural
Science Foundation. The work of Wang was also supported partly by the US National Science
Foundation (No.DMS-0072620).}

\author[DZC]{Daizhan Cheng},\quad
\author[DZC]{Lei Guo},\quad
\author[YL]{Yuandan Lin},\quad
\author[YL]{Yuan Wang}

\address[DZC]{Institute of  Systems Science,
Chinese Academy of Sciences, Beijing 100080, P.R.China\\
E-mail: dcheng@control.iss.ac.cn,\quad lguo@control.iss.ac.cn}

\address[YL]{Dept of Mathematical Sciences,
Florida Atlantic University, Boca Raton, FL 33431, USA\\
Email: lin@fau.edu, \quad ywang@fau.edu}

\begin{keyword}
Linear system, transition matrix, Squashing Lemma.
\end{keyword}

\begin{abstract}
In this note we show that for a given controllable pair $(A,B)$
and any $\lambda>0$, a gain matrix $K$ can be chosen so that the
transition matrix $e^{(A+BK)t}$ of the system $\dot x = (A+BK)x$
decays at the exponential rate $e^{-\lambda t}$ and the overshoot
of the transition matrix can be bounded by $M\lambda^L$ for some
constants $M$ and $L$ that are independent of $\lambda$.  As a
consequence, for any $h>0$, a gain matrix $K$ can be chosen so
that the magnitude of the transition matrix $e^{(A+BK)t}$ can be
reduced by $\frac{1}{2}$ (or by any given portion) over $[0, h]$.
An interesting application of the result is in the stabilization
of switched linear systems with any given switching rate (see
\cite{che04}).

\end{abstract}

\end{frontmatter}

\section{Introduction}

Consider a  linear system
\begin{equation}\label{e-sys}
\dot x = Ax+Bu,
\end{equation}
where $x(\cdot)$ takes values in $\R^n$,  $u(\cdot)$ takes values
in $\R^m$, and where $A$ and $B$ are matrices of appropriate
dimensions.  Suppose $(A, B)$ is a controllable pair. It is a well
known fact that for any $\lambda
>0$, a gain matrix $K$ can be chosen so that the transition matrix
of the system $\dot x = (A+BK)x$ decays exponentially at the rate
of $e^{-\lambda t}$, that is, for some $R>0$,
\[
\left\|e^{(A+BK)t}\right\| \le R e^{-\lambda t},
\]
where and hereafter $\|\cdot\|$ denotes the operator norm induced
by the Euclidean norm on $\R^n$. To get a faster decay rate, it is
natural to consider a ``higher gain'' matrix $K_1$. However, such
a gain matrix in general results in a bigger overshoot for the
transition matrix $e^{(A+BK_1)t}$. In this note, we show that in
the pole placement practice, a gain matrix $K$ can be chosen so
that the overshoot of the transition matrix $e^{(A+BK)t}$ can be
bounded by $M\lambda^L$ for some constants $M$ and $L$ independent
of $\lambda$.  As a consequence, one sees that for any $h>0$, a
gain matrix $K$ can be chosen so that the magnitude of the
transition matrix $e^{(A+BK)t}$ can be reduced by $\frac{1}{2}$
(or by any given portion) over $[0, h]$. Note that this is a
stronger requirement than merely requiring $e^{(A+BK)t}$ to decay
at an exponential rate.  An interesting application of the result
is in the stabilization of switched linear systems with a given
switching rate (see \cite{che04}).

The estimate of the overshoots of transition matrices in the
practice of pole assignments has been studied widely (see e.g.
\cite{loa77}, \cite{val01} and \cite{pai94}). Our main result in
this note can be considered an enhancement of the Squashing Lemma
(see \cite{pai94}, \cite{mor96} and \cite{hep99}) which says the
following: for any $\tau_0>0$, $\delta>0$, any $\lambda>0$, it is
possible to find $K$ such that
\begin{equation}\label{e-2}
\left\| e^{(A+BK)t}\right\|\le \delta e^{-\lambda(t-\tau_0)}.
\end{equation}
In the current note, we show that $K$ can be chosen so that the
estimate in (\ref{e-2}) can be strengthened to
\[
\left\| e^{(A+BK)t}\right\|\le M \lambda^Le^{-\lambda t}
\]
for some constants $M$ and $L$ which are independent of $\lambda$.
Our proof is constructive that shows explicitly how $M$ and $L$
are chosen.

\section{Main Result}

In this section we present our main result.

\noindent {\bf Proposition 2.1}\ Let $A\in \R^{n\times n}$ and
$B\in \R^{n\times m}$ be two matrices such that the pair $(A,B)$
is controllable. Then for any $\lmd > 0$, there exists a matrix $K\in
\R^{m\times n}$ such that 
\begin{align}
\label{3.1} \left\|e^{(A+BK)t}\right\|\leq M\lmd^L e^{-\lmd t},\quad
\forall\, t\geq 0,
\end{align}
where $L = (n-1)(n+2)/2$ and $M > 0$ is a constant, which is
independent of $\lmd$ and can be estimated precisely in terms of $A,
B$ and $n$.

Comparing with the Squashing Lemma obtained in \cite{pai94},
Proposition 2.1 has two improvements: (i). In \rref{e-2}, the
estimate on the transient overshoot is exponentially proportional
to the decay rate $\lambda$, which resulted in an estimation of
the transition matrix in terms of $e^{-\lambda(t-\tau_0)}$ instead
of $e^{-\lambda t}$. In \rref{3.1}, the estimate on the transient
overshoot is proportional to $\lambda^L$ instead of $e^{\lmd
\tau_0}$ as in \rref{e-2}. This distinction between the two types
of estimations may be significant for some possible extensions of
our results to systems with external inputs. (ii). The value of
the constant $M$ in estimate (\ref{3.1}) can be precisely
calculated by using our constructive proof (see equation
(\ref{3.131}) in the sequel). This is certainly a very desirable
feature for practical purposes. See Example 3.1 for some
illustrations.

Proposition 2.1 was primarily presented and applied to a
stabilization problem of switched linear systems in
\cite{cglw-ccc03}.  It was found later that a recent paper
\cite{fan02} also provides a similar result with similar proofs.
The difference is that \cite{fan02} only considered the single
input case and the upper bound $M\lmd^{L}$ in (\ref{3.1}) was found
to be a polynomial $p(\lmd)$ in \cite{fan02} without an explicit
expression.  Hence, our result has obvious merits in control
design.

\smallskip

\noindent {\bf Proof of Proposition 2.1.}\,\, First we consider a
linear system $(A, b)$ of a single input. Without loss of
generality, we assume that $(A, b)$ is in the Brunovsky canonical
form:
$$
A =\begin{pmatrix}
0&1 & 0 &\cdots & 0\\
0&0& 1 &\cdots & 0\\
\vdots & \vdots & \vdots &\ddots&\vdots\\
%{} & \cdots & {} &\cdots & 1\\
a_1 & a_2 &a_3&\cdots & a_n
\end{pmatrix}, \quad
b  = \begin{pmatrix}
  0\\ 0\\ \vdots\\
%  0\\
1
  \end{pmatrix}.
$$
Let $\lambda_1, \ldots, \lambda_n$ be $n$ distinct, negative real numbers.
There exists
some $k\in\R^{1\times n}$ such that the characteristic equation of the
closed-loop system $A + bk$ is $p(\lambda) = (\lambda-\lambda_1)(\lambda -
\lambda_2)\cdots (\lambda - \lambda_n)$.  Note that the closed-loop
system is given by
\begin{eqnarray*}
& &\dot x_1 = x_2, \quad\dot x_2 = x_3,\quad \ldots,\quad
\dot x_{n-1}=x_n,\\
& &\dot x_{n} = \beta_1 x_1 + \beta_2 x_2 + \cdots \beta_n x_n
\end{eqnarray*}
for some $\beta_1, \beta_2, \ldots, \beta_n\in\R$.
Hence, $x_1$ satisfies the equation
\begin{align}
\label{3.4} x_1^{(n)} = \beta_1 x_1 + \beta_2 \dot x_1 +\cdots +
\beta_n x_1^{(n-1)},
\end{align}
whose characteristic equation is the
%The corresponding characteristic equation of \rref{3.4} is the
same as $p(\lambda)$.  Hence, the general solution of
\rref{3.4} is
$$
x_1(t) = c_1e^{\lambda_1 t} + c_2 e^{\lambda_2t} + \cdots + c_n
e^{\lambda_n t},
$$
where $c_1, c_2, \ldots, c_n$ are constants. From the equations
$x_2 = \dot x_1, x_3 = \dot x_2, \ldots, x_n= \dot x_{n-1}$, we
have $x(t) = \Lambda_0 e^{Dt} c$, where
\begin{eqnarray*}
\Lambda_0  &=& \begin{pmatrix}
1&1&\cdots &1\\
\lambda_1&\lambda_2&\cdots & \lambda_n\\
\vdots &\vdots &\ddots&{}\\
\lambda_1^{n-1}&\lambda_2^{n-1}& \cdots &\lambda_n^{n-1}
  \end{pmatrix},\\
D  &=& \begin{pmatrix}
\lambda_1 & 0 &\cdots&0\\
0&\lambda_2 & \cdots & 0\\
\vdots & \vdots & \ddots &\vdots\\
0&0&\cdots&\lambda_n
  \end{pmatrix},
\end{eqnarray*}
and where $c  = \begin{pmatrix} c_1&c_2&\cdots c_n
  \end{pmatrix}^T$.
Now, observe that $x(0) = \Lambda_0 c$, that is, $c
=\Lambda_0^{-1}x(0)$ (note that $\Lambda_0$ is an invertible
Vandermonde matrix).  Comparing this with the transition matrix of
the system, one sees that
\begin{align}
\label{3.3} e^{(A + bk)t} = \Lambda_0 \,e^{Dt} \,\Lambda_0^{-1}.
\end{align}

Let $\lambda_{\max} = \max\{\abs{\lambda_1}, \ldots,
\abs{\lambda_n}\}$.  Without loss of generality, assume that
$\lambda_{\max}\ge 1$. To get an estimate on $\left\|\Lmd_0\right\|$
and $\left\| \Lmd_0^{-1}\right\|$,  we need the following simple
fact: for an $n\times n$ matrix $C$,  let $c_{\max}=\max_{1\leq
i,j\leq n}\abs{c_{ij}}$.  It is not hard to see that
\[
\norm{C}\le nc_{\max}.
\]
Hence, we have
\begin{align}
\label{3.6} \left\|\Lambda_0\right\|\leq n\lambda_{\max}^{n-1}.
\end{align}

To get an estimate on $\Lambda_0^{-1}$, first note that
\begin{equation}\label{e-L1}
\Lambda_0^{-1} = \frac{1}{{\rm det} \Lambda_0}\,  {\rm
adj}\,\Lambda_0,
\end{equation}
where ${\rm adj}\,\Lambda_0$ denotes the adjoint  matrix of
$\Lambda_0$, and that
$$
\displaystyle{\det \Lambda_0 = \prod_{j>i}(\lambda_j -
  \lambda_i)}.
$$
Hence, if we choose $\lambda_1, \ldots \lambda_n$ in such a way
that $\lambda_{i+1} \le \lambda_i - 1$  with $\lambda_1 <0$, we
get $\abs{\det \Lambda_0} \geq 1$.

Taking the structure of ${\rm adj \Lmd_0}$  into account, it is easy
to see that for $C={\rm adj}\Lmd_0$,
\begin{eqnarray*}
\phantom{xxxxx}c_{\max}&\leq& (n-1)!{\lmd_{\max}}^{1+2+\cdots+(n-1)}\\
&=& (n-1)!{\lmd_{\max}}^{\frac{n(n-1)}{2}}.
\end{eqnarray*}
Hence, by \rref{e-L1}, we have
\begin{align}
\label{3.8} \left\| \Lmd_0^{-1}\right\| \leq \norm{{\rm adj}\Lmd_0}
\leq n(n-1)!{\lmd_{\max}}^{\frac{n(n-1)}{2}}.
\end{align}
Consequently, \rref{3.6} and \rref{3.8} yield that
\begin{eqnarray*}
\left\| \Lambda_0\,e^{Dt}\,\Lambda_0^{-1}\right\|
%\left\|\Lambda_0\right\| \left\| e^{Dt}\right\|
%\left\|\Lambda_0^{-1}\right\| %\\\\ &\leq&
&\le& n \lambda_{\max}^{n-1}\,\left\| e^{Dt}\right\|\,n(n-1)!
\lambda_{\max}^{n(n-1)/2}\\
& \leq & nn!\lambda_{\max}^{(n-1)(n+2)/2}e^{-\lambda_{\min} t},
\end{eqnarray*}
where $\lambda_{\min} = \min\{\abs{\lambda_1}, \ldots,
\abs{\lambda_n}\}$.

Suppose for some $\rho >1$, $\lambda_{\max} \leq \rho
\lambda_{\min}$.  Then, it follows that
\begin{align}
\label{3.9}
 \left\|\Lambda_0e^{Dt}\Lambda_0^{-1}\right\|\leq M
\lambda_{\min}^{(n-1)(n+2)/2}e^{-\lambda_{\min} t},
\end{align}
where
 \begin{equation}\label{3.131}
 M=nn!\,\rho^{(n-1)(n+2)/2}.
 \end{equation}
In summary, we need the following conditions on the $\lambda_i$'s:

\begin{itemize}
\item $\lmd_1, \lmd_2, \cdots, \lmd_n$ are distinct, real, and negative;
\item $\lmd_{i+1}\leq \lmd_i - 1$ for $1\leq i\leq n-1$, and hence,
$\lmd_{\max} = \abs{\lmd_n}$, $\lmd_{\min} = \abs{\lmd_1}$;
\item
$\abs{\lmd_n}\leq \rho \abs{\lmd_1}$, for some constant $\rho
>1$.
\end{itemize}

Obviously, for any given $\lmd>0$, it is easy to choose $\lmd_1,
\cdots, \lmd_n$ to satisfy all the above conditions together with
the condition that $\lmd_1\le -\lmd$. For example, one can choose
$\lmd_1 <\min\{-1, -\lmd\}$,
 and let $\lmd_{i+1} = \lmd_i -1$ for $1\leq i \le n-1$. Since
$\abs{\lmd_n} = \abs{\lmd_1 -(n-1)}\le n \abs{\lmd_1}$, we see that
$\rho$ can be set as $\rho = n$.

With such choices of $\lambda_1, \lambda_2, \ldots, \lambda_n$, we
see from \rref{3.3} and \rref{3.9} that the desired result hold.

Now we consider the case when $(A, b)$ is not in the Brunovsky
canonical form.  In this case, find an invertible $T\in\R^{n\times
n}$ such that $(T^{-1}AT, T^{-1}b)$ is in the Brunovsky canonical
form.

For any given $\lmd>0$, the above proof has shown that for $A_1 =
T^{-1}AT$, $b_1 = T^{-1}b$, one can find $k_0\in \R^{1\times n}$
such that
\[
e^{(A_1+b_1k_0)t}\le M\lmd^Le^{-\lmd t},
\]
where $M$ is given by \rref{3.131} for some chosen $\rho$, and
$L=(n-1)(n+1)/2$.  Clearly, with $k=k_0T^{-1}$, one has
\begin{equation}\label{e-M1}
e^{(A+bk)t} = T(e^{(A_1+b_1k_0)t})T^{-1}\le M_1\lmd^Le^{-\lmd t},
\end{equation}
where $M_1 = M\norm{T}\norm{T^{-1}}$.

Finally, we consider the multi-input  system
\begin{align}
\label{3.12} \dot x = Ax + Bu,
\end{align}
where $A\in \R^{n\times n},  B \in \R^{n\times m}$.  Suppose that
the system is controllable.  By Heymann's Lemma (c.f., e.g., page
187 of \cite{son98}), one sees that for any $v\in\R^m$ such that
$b :=Bv\not=0$, there exists some $K_0\in\R^{m\times n}$ such that
$ (A + BK_0, \, b) $ is itself controllable.  Hence, the
conclusion of single-input case that has just been proved above is
applicable to the controllable pair $(A+BK_0, \, b)$, and one then
sees that there exists some $k\in \R^{1\times n}$ such that
$\left\|e^{(A + BK_0 + bk)t} \right\| \leq M\lmd^L e^{-\lmd t}$ for
all $t\ge 0$. Hence, with $K = K_0 + vk$, it holds that
\begin{align}
\label{3.13} \left\|e^{(A+BK)t}\right\|\leq M\lmd^L e^{-\lmd t}
\qquad\forall\,t\ge 0.
\end{align}
This completes the proof. \hfill $\Box$

\noindent{\bf Remark 2.2} In the above proof, we have used the
fact that for a single input system $(A, b)$ which is
controllable, when it is not in the Brunovsky canonical form, one
can find an invertible matrix $T$ such that $(T^{-1}AT, T^{-1}b)$
is in the canonical form.  To be more precise, the matrix $T$ can
be chosen as (see e.g., \cite{son98}):
\[
T=\begin{pmatrix}b&Ab&\cdots&A^{n-1}b\end{pmatrix}
\begin{pmatrix}a_{n-1}&\cdots& a_1 &1\\
\vdots&\vdots&1&0\\a_1&1&\cdots&0\\
1&0&\cdots&0 \end{pmatrix},
\]
where $a_1,\ldots, a_{n-1}$ are as in the characteristic
polynomial of $A$ given by
\[
 \det(sI-A) = s^n + a_1s^{n-1}+\cdots+a_{n-2}s+a_{n-1}.
\]
{}From this one can find an estimate of $\norm{T}$ and
$\norm{T^{-1}}$, which in turn will lead to an estimate of $M_1$
in \rref{e-M1}.

\section{An Example}

The design technique is demonstrated in the following example.

\noindent{\bf Example 3.1} Consider the following controllable
linear system:
$$
A_1=\begin{pmatrix} 1&0&1\\0&1&-1\\2&1&0
\end{pmatrix},\quad B_1=\begin{pmatrix} 1\\0\\1
\end{pmatrix},
$$
With the help of MATLAB, we first calculate the transfer matrix
$$
T_1 =\begin{pmatrix}0&0&1\\
-1&-1&0\\
-1&0&1
\end{pmatrix}.
$$
With the transfer matrix $T_1$, one has
$$
T_1^{-1}A_1T_1 =\begin{pmatrix}0&1&0\\
0&0&1\\
-1&0&2
\end{pmatrix},\quad
T_1^{-1}B_1 =\begin{pmatrix}0\\
0\\
1
\end{pmatrix}
$$
Calculation shows that $\|T_1\|=1.80193754431757$ and
$\|T_1^{-1}\|=2.24697960199992$. Taking $\rho=n(=3)$, we have
\begin{eqnarray}
& &L=\frac{(n-1)(n+2)}{2}=5,\label{e-new1}\\
& &M=\|T_1\|\|T_1^{-1}\|nn!n^{(n-1)(n-2)/2}\approx
218.642.\label{e-new2}
\end{eqnarray}

Suppose for some design purpose, a decay constant $\lmd = 49.894$ is
given. Choosing $\lmd_1=-\lmd$, $\lmd_2=\lmd_1-1$, $\lmd_3=\lmd_2-1$, the
feedback $K_1$ can be easily calculated (under the normal form) as
$$
\tilde{K}_1\approx \begin{pmatrix} -151.681 &-7769.474
&-131773.562\end{pmatrix}.
$$
Back to the original coordinate frame, we have
$$
K_1=\tilde{K}_1T_1^{-1}\approx\begin{pmatrix}
  -124155.769 & 7769.474&-7617.793
\end{pmatrix}.
$$
With such a choice of $K_1$, we get the desired decay estimate
\[
\norm{e^{(A+BK)t}}\le M\lmd^Le^{-\lmd t}\quad\forall\,t\ge 0,
\]
for the given decay constant $\lmd = 49.894$ with $L$ and $M$ given
as in \rref{e-new1}--\rref{e-new2}. \hfill $\Box$

\comment{
\noindent{\bf Remark 3.2.} In general, the amplitude of
the control law $K_1$ may become large when $\lmd_1$ is large. This
seems to be inevitable in practice since high gain controller has
significant advantages in dealing with uncertainties in the system
structure.
}

\section{Conclusion}

In this note we show that if $(A, B)$ is controllable, then for
any $\lmd>0$, a gain matrix $K$ can be chosen such that the
transition matrix $e^{(A+BK)t}$ decays at the exponential rate
$e^{-\lmd t}$ and the overshoot of $e^{(A+BK)t}$ can be bounded by
$M\lmd^L$ for some constants $M$ and $L$ that are independent of the
decay constant $\lmd$.  The result provides a convenient tool for
control design, particularly for switched systems, see
\cite{che04}.

\vspace{0.4in}

\centerline{\bf Erratum}

\medskip

There is a mild flaw in the statement of Proposition 2.1 in the above paper
(cf.~\cite{sqlemma1}). We restate it as follows.

\noindent{\bf Proposition}\
\ Let $A\in \R^{n\times n}$ and
$B\in \R^{n\times m}$ be two matrices such that the pair $(A,B)$
is controllable. Then for any $\lmd
\ge 1$, there exists a matrix $K\in \R^{m\times n}$ such that
\begin{align}
%\label{3.1} 
\left\|e^{(A+BK)t}\right\|\leq M\lmd^L e^{-\lmd t},\quad
\forall\, t\geq 0,
\end{align}
where $L = (n-1)(n+2)/2$ and $M > 0$ is a constant, which is
independent of $\lmd$ and can be estimated precisely in terms of $A,
B$ and $n$.~\mbox

The proof of Proposition 2.1 in \cite{sqlemma1} is only valid for the
case when $\lambda \ge 1$ (instead of the original version of
$\lambda>0$), because the eigenvalues $\lambda_1,\ldots, \lambda_n$ were
chosen to satisfy $\lambda_1 \le -1$, and $\lambda_k\le\lambda_1$ for
$k\ge 1$.  For more details, we refer the reader to the discussions that
followed formula (10) in \cite{sqlemma1}.

A main motivation of the work in \cite{sqlemma1} was for us to develop
the results in \cite{switching}. As in most applications of overshoot
estimation for pole placements, the parameter $\lambda$ in
\cite{switching} was chosen as a number of large value.  Hence, the
correction does not affect our results in \cite{switching}.

\noindent
{\it Acknowledgment.}  The authors would like to thank Prof.~Elena De
Santis for pointing out the error.

\end{document}